\documentclass[conference]{IEEEtran}
\IEEEoverridecommandlockouts
\usepackage{cite}
\usepackage{amsmath,amssymb,amsfonts}
\usepackage{algorithmic}
\usepackage{graphicx}
\usepackage{textcomp}
\usepackage{xcolor}
\usepackage{hyperref}
\hypersetup{colorlinks=true,
			urlcolor=blue,
			linkcolor=black,
			citecolor=blue,
			bookmarksdepth=paragraph}

\DeclareMathOperator*{\argmin}{arg\,min}

\def\BibTeX{{\rm B\kern-.05em{\sc i\kern-.025em b}\kern-.08em
    T\kern-.1667em\lower.7ex\hbox{E}\kern-.125emX}}
\begin{document}

\title{Synthesizing Control Laws from Data using Sum-of-Squares Optimization\\
\thanks{This work was supported by the NSERC discovery grant program.}
}

\author{\IEEEauthorblockN{Jason J. Bramburger}
\IEEEauthorblockA{\textit{Dept.~of Mathematics and Statistics} \\
\textit{Concordia University}\\
Montr\'eal, Canada \\
jason.bramburger@concordia.ca}
\and
\IEEEauthorblockN{Steven Dahdah}
\IEEEauthorblockA{\textit{Dept.~of Mechanical Engineering} \\
\textit{McGill University}\\
Montr\'eal, Canada \\
steven.dahdah@mail.mcgill.ca}
\and
\IEEEauthorblockN{James Richard Forbes}
\IEEEauthorblockA{\textit{Dept.~of Mechanical Engineering} \\
\textit{McGill University}\\
Montr\'eal, Canada \\
james.richard.forbes@mcgill.ca}}

\maketitle

\begin{abstract}
The control Lyapunov function (CLF) approach to nonlinear control design is well established. Moreover, when the plant is control affine and polynomial, sum-of-squares (SOS) optimization can be used to find a polynomial controller as a solution to a semidefinite program. This letter considers the use of data-driven methods to design a polynomial controller by leveraging Koopman operator theory, CLFs, and SOS optimization. First, Extended Dynamic Mode Decomposition (EDMD) is used to approximate the Lie derivative of a given CLF candidate with polynomial lifting functions. Then, the polynomial Koopman model of the Lie derivative is used to synthesize a polynomial controller via SOS optimization. The result is a flexible data-driven method that skips the intermediary process of system identification and can be applied widely to control problems. The proposed approach is used to successfully synthesize a controller to stabilize an inverted pendulum on a cart.
\end{abstract}



\section{Introduction}

Nonlinear systems can be found throughout engineering, science, economics, and other domains. Often nonlinear systems must be controlled to realize useful behavior. There are a plethora of nonlinear control design methods to choose from, such as gain-scheduling, feedback linearization, integrator backstepping, sliding-mode control, and others \cite{book_khalil}. 

A popular nonlinear control design method is the \emph{control Lyapunov function (CLF)} approach \cite{sontag2013mathematical}. The gist of the CLF approach to control design is that, given a control affine system, a controller is sought such that a Lyapunov function candidate associated with the closed-loop system satisfies $V(x_*) = 0$, $V(x) > 0, \, \forall x \in \mathcal{D} \setminus \left\{x_* \right\}$, and $\dot{V}(x) < 0, \, \forall x \in \mathcal{D} \setminus \left\{x_* \right\}$ where $x_*$ is an equilibrium point and $\mathcal{D} \subseteq \mathbb{R}^d$ is some domain. See \S\ref{sec:cont_Lyap_fun} for a review of the CLF approach to controller design.

One means of designing a controller via the CLF approach is to employ techniques from \emph{sum-of-squares (SOS) optimization}. Roughly speaking, by exploiting the polynomial form of the plant and the Lyapunov function candidate, a sufficient condition for the existence of a polynomial controller is the existence of a solution to a linear matrix inequality (LMI) feasibility problem \cite{anderson2015advances,papachristodoulou2005tutorial}. The SOS optimization approach to control design via a CLF is attractive because semidefinite programming can be leveraged to solve LMI feasibility problems in a simple and efficient manner.

There have been numerous SOS optimization approaches to CLF-based control design. For instance, \cite{Wloszek_2003,Majumdar_2013} consider controller design using a CLF approach to ensure that the region of attraction of the closed-loop system about the equilibrium point $x_*$ is as large as possible. In \cite{yang2023suboptimal}, an approximate solution to the Hamilton-Jacobi-Bellman (HJB) equation is found using SOS optimization, yielding a suboptimal controller. State- and output-feedback control design in a SOS optimization framework for parabolic PDE systems is considered in \cite{Gahlawat_2017}. 

To use the CLF approach in concert with SOS optimization for control design, a model of the nonlinear system is needed. When a model is not available, but a plethora of data is available, a data-driven approach to modelling and control design is natural. The \emph{Koopman operator} approach to data-driven modelling of nonlinear systems has garnered significant attention recently \cite{mauroy_2020_koopman,otto_2021_koopman,brunton_2021_modern}. The basic idea behind the Koopman operator is that a finite-dimensional nonlinear system can be expressed as an infinite-dimensional linear system using \textit{lifting functions} \cite{koopman_hamiltonian_1931}. The linearity of the Koopman operator is attractive because standard linear systems tools, such as the eigenspectrum \cite{mezic_2019_spectrum}, can be used to analyze nonlinear systems. A finite-dimensional approximation of the Koopman operator can be readily identified from data~\cite{eDMD}. This approximate Koopman operator can then be used as the basis for control design. For instance, \cite{korda_2018_linear} considers model predictive control (MPC),  \cite{abraham_active_2019} considers an active learning approach in a Koopman framework, and \cite{moyalan_2022_data-driven} considers discounted optimal control using the Perron-Frobenius operator, the dual to the Koopman operator.

This paper proposes a data-driven approach to CLF-based control design using the Koopman operator. The CLF approach relies on the generator of the Koopman operator, the so-called Lie derivative, of the closed-loop Lyapunov function candidate being strictly negative. When the lifting functions associated with the Koopman operator are polynomials, the Lie derivative of the closed-loop Lyapunov function candidate is \emph{almost} polynomial. Forcing the controller to be polynomial enables the use of SOS optimization to find a suitable controller that renders the Lie derivative of the closed-loop Lyapunov function candidate negative definite. As a consequence, the equilibrium point $x_*$ of the closed-loop system is made asymptotically stable.

The novel contribution of this letter is nonlinear controller synthesis via the following two-step method. First, the Lie derivative of a given CLF candidate is approximated using the Koopman operator with polynomial lifting functions. Next, this Koopman representation of the Lie derivative is incorporated into a SOS optimization problem that parameterizes the controller as a polynomial of the state variables. Convex SOS optimization routines are then leveraged to find a control law that renders the equilibrium point $x_*$ of the closed-loop system asymptotically stable.

This letter is organized as follows. CLFs, SOS optimization, and Koopman operator theory are reviewed in \S\ref{sec:preliminaries}. The main theoretical results are presented in \S\ref{sec:main_results}. A numerical example involving the control of an inverted-pendulum on a cart is provided in \S\ref{sec:inv_pend_ex}. The paper is drawn to a close in \S\ref{sec:conclusion}.


\section{Preliminaries}
\label{sec:preliminaries}

In this section we provide the necessary preliminary information to present our method of synthesizing control laws from data. Throughout we will have $\langle u,v\rangle$ denote the inner product of vectors $u,v\in\mathbb{R}^d$. The temporal argument of the state $x(t)$, the control $u(t)$, etc., will be suppressed unless required for clarity.

\subsection{Control Lyapunov functions}
\label{sec:cont_Lyap_fun}

Consider a control affine system
\begin{equation}\label{ControlAffine}
	\dot{x} = f(x) + \sum_{i = 1}^m g_i(x)u_i, \quad x \in \mathbb{R}^d, \ u_i \in \mathbb{R},
\end{equation}
for which we wish to specify a control input $u = [u_1, \dots, u_m]^T$ that forces all initial conditions belonging to some domain $\mathcal{D} \subseteq \mathbb{R}^d$ into an equilibrium point $x_*$ as $t \to \infty$ under the flow of \eqref{ControlAffine}. To do so, one may specify a {\em control Lyapunov function} (CLF) \cite{sontag2013mathematical} $V: \mathcal{D} \to \mathbb{R}$ satisfying $V(x) > 0 \, \forall x \in \mathcal{D}\setminus\{x_*\}$ and $V(x_*) = 0$, and then seek a control input so that 
\begin{equation}\label{NegativeLyap}
    \bigg\langle \nabla V(x),f(x) + \sum_{i = 1}^m g_i(x)u_i\bigg\rangle < 0, \quad \forall x\in\mathcal{D}\setminus\{x_*\}.
\end{equation}
Indeed, \eqref{NegativeLyap} guarantees that $V$ decreases monotonically along trajectories of \eqref{ControlAffine}, eventually reaching the global minimum at $V(x_*) = 0$. Crucial for our work later in this letter, the system being control affine means that $u$ enters \eqref{NegativeLyap} linearly. 

When the inequality in \eqref{NegativeLyap} is not strict, only Lyapunov stability can be concluded rather than asymptotic stability. Precisely, this means that $0 \leq V(x(t)) \leq V(x(0))$ for all $t \geq 0$, meaning that the motion of $x(t)$ is constrained by the level set $V(x(0))$. Moreover, LaSalle's invariance principle guarantees the limiting behaviour of $x(t)$ is contained in the set of $x$ values for which \eqref{NegativeLyap} is exactly zero. Such non-strict inequalities become important when imposing a tractable relaxation of the inequality \eqref{NegativeLyap} in the following subsection.

\subsection{Synthesizing controllers with semidefinite programming}
\label{sec:SDPs}

Although finding a control input $u$ to satisfy \eqref{NegativeLyap} is difficult in general, under certain assumptions on \eqref{ControlAffine} this process can be automated using standard optimization methods. Specifically, if $f(x)$, $g_i(x)$, and $V(x)$ are all polynomial, then we can search for a polynomial state-dependent feedback control law $u(x)$. Fixing the degree of $u(x)$ allows its coefficients to be optimized to satisfy the polynomial inequality constraint \eqref{NegativeLyap} for all $x \in \mathcal{D}$. Although these coefficients appear linearly in \eqref{NegativeLyap}, tuning them to satisfy such a polynomial inequality is an NP-hard task in general \cite{murty1987some}.  

To make the controller synthesis problem tractable, we replace the polynomial inequalities with the sufficient conditions that the polynomials are sum-of-squares (SOS). Precisely, a polynomial $p(x)$ is SOS if there exists polynomials $q_1(x),\dots,q_k(x)$ so that 
\begin{equation}
    p(x) = \sum_{i = 1}^k q_i(x)^2.    
\end{equation}
Identifying an SOS representation of a polynomial trivially verifies that it is nonnegative. Furthermore, verifying an SOS representation constitutes a semidefinite program, since $p(x)$ is SOS if and only if there exists a vector of monomials $v(x)$ and a positive semidefinite matrix $P$ such that $p(x) = v(x)^TPv(x)$.

Recall \eqref{NegativeLyap} and suppose that $\mathcal{D}$ is a semialgebraic set, meaning there exists polynomials $\{a_j\}_{j = 1}^J$, $\{b_\ell\}_{\ell = 1}^L$ so that
\begin{equation}\label{Semialg}
    \mathcal{D} = \{x \in \mathbb{R}^d|\ a_j(x) \geq 0, b_\ell(x) = 0, \forall j,\ell\}.
\end{equation}
A sufficient condition for identifying a stabilizing polynomial state-dependent controller $u(x) = \begin{bmatrix} u_1(x), \dots, u_m(x) \end{bmatrix}^T$ is~\cite{anderson2015advances,papachristodoulou2005tutorial}
\begin{equation}\label{SOSLyap}
    \begin{split}
        &(a)\ -\bigg\langle \nabla V(x),f(x) + \sum_{i = 1}^m g_i(x)u_i(x)\bigg\rangle \\ 
        &\qquad + \sum_{j = 1}^Ja_j(x)\sigma_j(x) + \sum_{\ell = 1}^L b_\ell(x)\rho_\ell(x)\ \mathrm{is\ SOS},\\
        &(b)\ \sigma_j(x)\ \mathrm{is\ SOS}, 
    \end{split}
\end{equation}
for polynomials $\sigma_j$ and $\rho_\ell$. Note that ``is SOS'' means there exists an SOS representation of the given polynomial. By fixing the degrees of $u_i$, $\sigma_j$, and $\rho_\ell$, the SOS constraints can be translated into semidefinite programs by freely available software packages like YALMIP~\cite{lofberg2004yalmip}. Solvers such as MOSEK~\cite{mosek} can efficiently solve these semidefinite programs to determine the polynomial coefficients.

\subsection{The Koopman operator}
\label{sec:koopman_review}

Let $\Phi(t;x):\mathbb{R}_+\times X \to X$ represent the flow of a dynamical system at time $t \geq 0$ with the initial condition $\Phi(0;x) = x$. The {\em Koopman operator} is a linear transformation that $\mathcal{K}_t$ that for each $t$ maps a lifting function $\varphi:X \to \mathbb{R}$ to
\begin{equation}\label{Koopman}
    \mathcal{K}_t\varphi = \varphi(\Phi(t;x)), \quad \forall x \in \mathbb{R}^d.
\end{equation}
Koopman lifting functions are also often called \textit{observables}, but are unrelated to the concept of observability.
The family $\{\mathcal{K}_t|\ t \in \mathbb{R}_+\}$ is a one-parameter semigroup whose generator is referred to as the \textit{Lie derivative}, acting on differentiable lifting functions via
\begin{equation}\label{Lie}
    \mathcal{L}\varphi = \lim_{t \to 0^+} \frac{\mathcal{K}_t \varphi - \varphi}{t}.
\end{equation}
Intuitively, the Lie derivative represents a derivative of the lifting function $\varphi$ in the direction of the flow of the system. 

Lyapunov functions are examples of lifting functions whose values decrease along trajectories. Moreover, given a candidate Lyapunov function $V:\mathcal{D}\to\mathbb{R}$ for \eqref{ControlAffine}, the Lie derivative of $V$ is evaluated using the chain rule to be exactly
\begin{equation}\label{controlLie}
    \mathcal{L}V = \bigg\langle \nabla V(x),f(x) + \sum_{i = 1}^m g_i(x)u_i\bigg\rangle, 
\end{equation}
which demonstrates a connection between the Koopman operator and Lyapunov functions.

\section{Main Result}
\label{sec:main_results}

Our method of synthesizing control laws from data comes as a two-step process. First, we estimate the Lie derivative from data using well-developed techniques for approximating the action of the Koopman operator on finite \textit{dictionaries} of lifting functions \cite{eDMD}. Second, we integrate our estimated Lie derivative into the SOS framework of \S\ref{sec:SDPs} and provide an SOS optimization problem that determines control laws directly from data.

\subsection{Estimating Lie derivatives from data}

Begin by supposing that snapshots of a control system are given in the form of triples $\{(x_k,u_k,y_k)\}_{k = 1}^n \subset \mathbb{R}^d\times\mathbb{R}^m\times\mathbb{R}^d$. Here, $y_k$ is the state of the system exactly $\tau > 0$ time units after having state and control values $(x_k,u_k)$. We then consider two dictionaries of polynomial lifting functions of the state, $\phi_1,\dots,\phi_p$ and $\psi_1,\dots,\psi_q$. Denote
\begin{equation}
    \phi = \begin{bmatrix} \phi_1 & \cdots & \phi_p \end{bmatrix}^T, \quad \psi = \begin{bmatrix} \psi_1 & \cdots & \psi_q \end{bmatrix}^T .
\end{equation}  
Lyapunov functions are assumed to belong to $\mathrm{span}\{\phi\}$, while their image under the Lie derivative is projected into $\mathrm{span}\{\psi\}$. Two dictionaries are required because the Lie derivative associated with a polynomial system is expected to be of a higher degree than the original lifting function it is being applied to. To see why this is the case, let $F$ denote the right-hand-side of our (assumed polynomial) control-affine system \eqref{ControlAffine}. Note that the degree of $\mathcal{L}V = \langle \nabla V,F\rangle$, the Lie derivative \eqref{controlLie}, is 
\begin{equation}
    \mathrm{deg}(\nabla V) + \mathrm{deg}(F) = \mathrm{deg}(V) - 1 + \mathrm{deg}(F) \geq \mathrm{deg}(V),   
\end{equation}
with a strict inequality when $F$ is nonlinear, that is, $\mathrm{deg}(F) > 1$. In practice one should aim to have $\mathrm{span}\{\phi\} \subseteq \mathrm{span}\{\psi\}$ for the approximation of the Lie derivative later in this section.   

Following the method of EDMD \cite{eDMD}, define the $p\times n$ matrix
\begin{equation}
	\Phi = \begin{bmatrix}
		\phi(y_1) & \phi(y_2) & \cdots & \phi(y_n)
	\end{bmatrix},
\end{equation}
and the $q (m+1)\times n$ matrix
\begin{equation}
	\Psi = \begin{bmatrix}
		\psi(x_1) & \psi(x_2) & \cdots & \psi(x_n) \\
		\psi(x_1)u_{1,1} & \psi(x_2)u_{1,2} & \cdots & \psi(x_n)u_{1,n} \\
        \vdots & \vdots & \ddots & \vdots \\
		\psi(x_1)u_{m,1} & \psi(x_2)u_{m,2} & \cdots & \psi(x_n)u_{m,n} 
	\end{bmatrix},
\end{equation}
where we assume that the dynamics are control affine, leading to the specific form of $\Psi$~\cite{bruder_2021_advantages}. With these matrices we can approximate the Koopman operator by first obtaining the matrix $K \in \mathbb{R}^{p \times q(m+1)}$ as the solution to the minimization problem 
\begin{equation}
    K = \Phi\Psi^\dagger = \argmin_K \|\Phi - K\Psi\|_F,
\end{equation}
where ${(\cdot)}^\dagger$ denotes the Moore-Penrose pseudoinverse and $\|\cdot\|_F$ denotes the Frobenius norm. The matrix $K$ can be broken up into $m+1$ size $p \times q$ matrices 
\begin{equation}\label{KoopData}
    K = \begin{bmatrix} A & B_1 & \cdots & B_m \end{bmatrix},
\end{equation}
representing the lifted components corresponding to $f$ and $g_1,\dots,g_m$ in a control affine system \eqref{ControlAffine}. Then, for any given control input $u = \begin{bmatrix} u_1, \dots, u_m \end{bmatrix}^T$, the matrix $K$ leads to an approximation of the Koopman operator $\tilde{\mathcal{K}}$ acting on lifting functions $\varphi = \langle c, \phi\rangle \in \mathrm{span}\{\phi\}$ by
\begin{equation}\label{KoopApprox}
    \tilde{\mathcal{K}}(\varphi) := \langle c, A\psi\rangle + \sum_{i = 1}^m\langle c, B_i\psi\rangle u_i, 
\end{equation}
which one can verify belongs to $\mathrm{span}\{\psi\}$ for each $u \in \mathbb{R}^m$. 

With the approximate Koopman operator $\tilde{\mathcal{K}}$ we can further estimate the Lie derivative as a finite-difference approximation of \eqref{Lie}. Precisely, $\tilde{\mathcal{L}}$ acts on $\varphi = \langle c, \phi\rangle \in \mathrm{span}\{\phi\}$ by 
\begin{equation}\label{LieData}
    \tilde{\mathcal{L}}(\varphi) = \frac{\tilde{\mathcal{K}}(\varphi) - \varphi}{\tau}.
\end{equation}
Then, using \eqref{KoopApprox}, it follows that
\begin{equation}
    \tilde{\mathcal{L}}(\varphi) = \tau^{-1}\langle c, A\psi - \phi\rangle + \tau^{-1}\sum_{i = 1}^m\langle c, B_i\psi\rangle u_i,        
\end{equation}
for a given control input $u \in \mathbb{R}^m$. Having $\mathrm{span}\{\phi\} \subseteq \mathrm{span}\{\psi\}$ guarantees that $\tilde{\mathcal{L}}:\mathrm{span}\{\phi\} \to \mathrm{span}\{\psi\}$, just like $\tilde{\mathcal{K}}$.  We refer the reader to \cite{bramburger2023auxiliary} for convergence proofs regarding the Lie derivative approximation $\tilde{\mathcal{L}}$ in the limits of infinite data ($n\to \infty$), sampling rates ($\tau \to 0^+$), and dictionaries ($q \to \infty$).

\subsection{Synthesizing control laws from data}

The work in the previous subsection allows one to estimate the Lie derivative direction from data. We now leverage this method to incorporate it into a SOS-driven method for synthesizing control laws. Let us begin by assuming that $V \in \mathrm{span}\{\phi\}$, that is, there is some $c \in \mathbb{R}^p$ so that $V = \langle c,\phi\rangle$ is our candidate control Lyapunov function. The Lie derivative condition \eqref{NegativeLyap} is then replaced with the data-driven Lie derivative condition,
\begin{equation}
    \tilde{\mathcal{L}}(V) = \tau^{-1}\langle c, A\psi - \phi\rangle + \tau^{-1}\sum_{i = 1}^m\langle c, B_i\psi\rangle u_i < 0.
    \label{eq:data-driven_Lie_deriv_cond}
\end{equation}
Since $c$ is fixed because $V$ is given, \eqref{eq:data-driven_Lie_deriv_cond} is an affine constraint for the control input $u$. Moreover, since $\phi$ and $\psi$ are dictionaries of polynomial lifting functions, it follows that all terms in $\tilde{\mathcal{L}}(V)$ are polynomials in the state variable $x$. Thus, we now follow a similar procedure to \S\ref{sec:SDPs} by considering $u$ as a polynomial function of the state variable $x$ and relaxing the inequalities to SOS conditions.

In detail, we consider a third dictionary of polynomials in $x$, denoted
\begin{equation}
    \chi = \begin{bmatrix} \chi_1,\dots,\chi_r \end{bmatrix}^T,
\end{equation}
and consider $u = C\chi$ for some coefficient matrix $C \in \mathbb{R}^{m \times r}$. Hence, the data-driven SOS relaxation for synthesizing control laws on the semialgebraic set $\mathcal{D}$ as in \eqref{Semialg} is given by
\begin{equation}\label{SOSLyapData}
    \begin{split}
        &(a)\ -\tau^{-1}\langle c, A\psi - \phi\rangle - \tau^{-1}\sum_{i = 1}^m\langle c, B_i\psi\rangle [C\chi]_i  \\ 
        &\qquad + \sum_{j = 1}^Ja_j(x)\sigma_j(x) + \sum_{\ell = 1}^L b_\ell(x)\rho_\ell(x)\ \mathrm{is\ SOS},\\
        &(b)\ \sigma_j(x)\ \mathrm{is\ SOS}, 
    \end{split}
\end{equation}
which comes from replacing the exact Lie derivative, $\mathcal{L}V$ in \eqref{SOSLyap}, with that approximated from data, $\tilde{\mathcal{L}}(V)$, in \eqref{eq:data-driven_Lie_deriv_cond}. The goal is then to determine the coefficients of the matrix $C$ appropriately to derive a control law from data. Since there could be many choices for $C$, we propose the following convex optimization problem:
\begin{equation}\label{ControlMin}
    \min_{C \in\mathbb{R}^{m \times r}} \{h(C): \eqref{SOSLyapData}\ \mathrm{is\ satisfied}\}
\end{equation}
where $h:\mathbb{R}^{m\times r} \to \mathbb{R}$ is a convex optimization objective.

The optimization objective is user-specified and should be guided by the specifics of the problem. Well-motivated examples of $h:\mathbb{R}^{m\times r} \to \mathbb{R}$ are as follows. 
\begin{itemize}
    \item \textit{Sparsity} --- A convex proxy for producing control laws that use the fewest elements of $\chi$ possible can be implemented by setting $h(C) = \|C\|_1$, the absolute sum of all elements in $C$. The resulting controller relies on the fewest elements of $\chi$ and may help to tame the influence of noisy measurements. 
    \item \textit{Boundedness} --- Physical limitations of actuators often require that the controller output cannot exceed certain values. Synthesizing a controller that has bounded fluctuation over the set $\mathcal{D}$ can be implemented by minimizing $h(C) = M$ so that $M - u_i \geq 0$ and $u_i - M \geq 0$ for all $x \in \mathcal{D}$. Indeed, these conditions guarantee that $|u_1|,\dots,|u_m| \leq M$ for all $x \in \mathcal{D}$, while the inequalities can be relaxed to SOS conditions on $\mathcal{D}$ similar to \eqref{SOSLyap} and appended to the conditions in \eqref{ControlMin}. One can similarly impose boundedness on derivatives of $u$ in the same way, thereby enforcing a bound on the rate of control.
\end{itemize}


\section{Application to the inverted pendulum}
\label{sec:inv_pend_ex}

As a demonstration of the efficacy of our proposed method, we apply it to synthetic data generated from a simple inverted pendulum on a cart. Let $\theta$ denote the angle of the pendulum arm, with $\theta = 0$ denoting the upright position. The motion of the pendulum arm can be derived from first principles and is captured by the control affine system \cite{maeba2010swing}
\begin{equation}\label{PendCart}
	\ddot\theta = \sin(\theta) - \varepsilon \dot\theta - \cos(\theta)u,
\end{equation}
where $\varepsilon > 0$ is the scaled viscous friction coefficient and the control input $u$ is the acceleration of the cart on which the pendulum is placed. We fix $\varepsilon = 0.1$ throughout our demonstrations. Motivated by the control Lyapunov functions presented in \cite{maeba2010swing}, we choose our candidate Lyapunov function to be
\begin{equation}\label{PendLyap}
	V(\theta,\dot\theta) = \frac{1}{2}\dot\theta^2 + 1 - \cos(\theta) + \alpha(1 - \cos^3(\theta)), 
\end{equation} 
for some $\alpha > 0$. This is, of course, only one of many possible Lyapunov functions, and experiments with other Lyapunov functions are similarly promising. For the purposes of reproducing our results, all code related to this demonstration can be found at \href{https://github.com/jbramburger/data-clf}{github.com/jbramburger/data-clf}, which contains both MATLAB and Python implementations of the method. 

Synthetic data is produced by simulating 20 random initial conditions $(\theta(0),\dot{\theta}(0)) \in [0,2\pi)\times[-2,2]$ subject to a sinusoidal external forcing of the form
\begin{equation}
	u(t) = A\sin(t + B), 
\end{equation} 
where $A$ and $B$ are drawn from the uniform distribution on $[-1,1]$ and $[-\pi,\pi]$, respectively. We integrate each initial condition from $t = 0$ up to $t = 20$, collecting data at evenly spaced intervals of length $\tau = 0.01$. It should be noted that with a physical pendulum on a cart one can easily perform such experiments quickly and collect the resulting data for $\theta(t)$, $\dot\theta(t)$, and $u(t)$ \cite{kaheman2022experimental}. Collecting data directly from the experiment could also improve control results as models such as \eqref{PendCart} are only approximate, thus circumventing the need for parameter estimation from a physical model.  

Considering the phase variable $\theta$ modulo $2\pi$ constrains the dynamics of the pendulum to a cylinder parameterized by the state variables $(\theta,\dot\theta)$. We embed this cylinder in $\mathbb{R}^3$ by introducing the lifted state variables
\begin{equation}
    (x_1,x_2,x_3) = (\cos(\theta),\sin(\theta),\dot\theta).
\end{equation}
In these lifted variables the Lyapunov function \eqref{PendLyap} becomes
\begin{equation}
	V(x_1,x_2,x_3) = \frac{1}{2}x_3^2 + 1 - x_1 + \alpha(1 - x_1^3),
\end{equation}
which is now a polynomial and has a global minimum at the unstable equilibrium $x_* = (1,0,0)$. Due to the restriction $1 - x_1^2 - x_2^2 = 0$, our dictionaries need only have linear terms in $x_2$ since higher-order terms can equivalently be provided through constant terms and powers of $x_1$. With this in mind, we take $\phi$ to be a dictionary of all monomials $x_1^{\alpha_1}x_2^{\alpha_2}x_3^{\alpha_3}$ with $\alpha_1,\alpha_3 = 0,1,2,3$ and $\alpha_2 = 0,1$, so that $V \in \mathrm{span}\{\phi\}$. For demonstration, we take $\psi$ to be a dictionary of monomials for which $\alpha_1,\alpha_3 = 0,1,2,3,4$ and $\alpha_2 = 0,1$, although using larger values of $\alpha_1,\alpha_3$ returns similar results.

To implement the synthesis of a control law $u(x_1,x_2,x_3)$ as an optimization problem, we use the state space 
\begin{equation}
    \mathcal{D} = \{(x_1,x_2,x_3)\in\mathbb{R}^3|\ \eta^2 - x_2^2 \geq 0, \ 1 - x_1^2 - x_2^2 = 0\}
\end{equation}
for some parameter $\eta \in (0,1)$, which must be built into our SOS programs as in \eqref{SOSLyap} and \eqref{SOSLyapData}. 
The equality condition $1 - x_1^2 - x_2^2 = 0$ was discussed above, while the inequality condition $\eta^2 - x_2^2 \geq 0$ excludes from the domain a strip on the cylinder where $x_1 = 0$ corresponding to $\theta = \frac{\pi}{2},\frac{3\pi}{2}$. The reason that this is excluded from the domain is that the control input to the system, given by $\cos(\theta)u = x_1u$, vanishes at this point and the Lie derivative is not necessarily negative for all $x_2 = \pm 1$ and $x_3 \in \mathbb{R}$. In practice, we take $\eta$ as close to 1 as possible to encapsulate as much of the the full state space as we can. Numerical results presented in this demonstration use $\eta^2 = 0.95$ and we promote sparsity in the controller using the optimization objective $h(C) = \|C\|_1$, as presented previously. 

\begin{figure}[htbp] 
\center
\includegraphics[width = 0.48\textwidth]{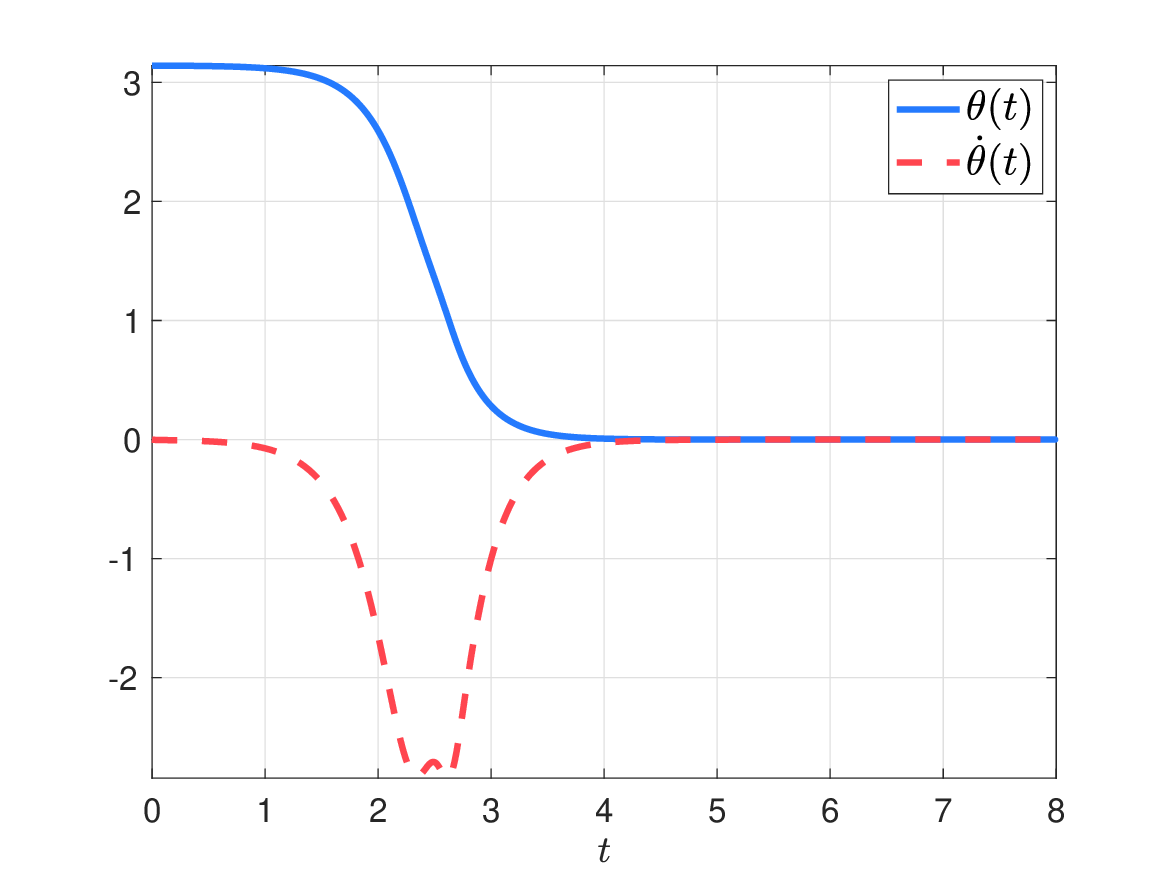} \\ 
\includegraphics[width = 0.48\textwidth]{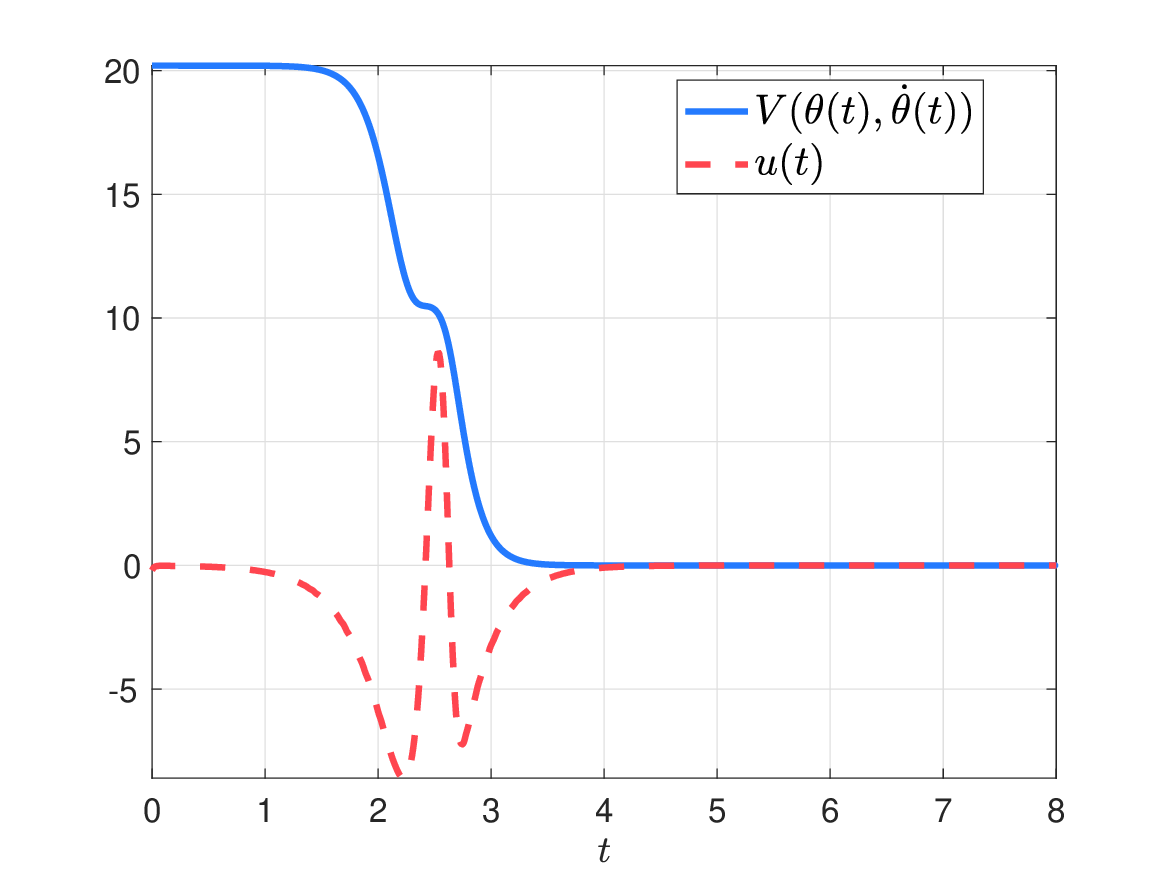}  
\caption{Controlled solutions of the inverted pendulum on a cart model \eqref{PendCart}. Top: The controlled arm angle $\theta$ and its derivative $\dot\theta$, forced to converge to the unstable upright position at $(\theta,\dot\theta) = (0,0)$. Bottom: The monotonic approach inside the Lyapunov function \eqref{PendLyap} (scaled by $1/10$ for interpretability) according to the controller, $u(\theta,\dot\theta)$, given by \eqref{PendController}.}
\label{fig:PendControl}
\end{figure} 

Experiments reveal that increasing $\alpha$ causes a proportional increase in the exponential rate of convergence of the pendulum arm to the upright position $(\theta,\dot\theta) = (0,0)$. As an example, Figure~\ref{fig:PendControl} presents controlled solutions with $\alpha = 100$. The initial condition is taken to be close to the hanging down position so that the movement of the cart is forced to swing the pendulum up into the upright position. The synthesized state-dependent control law is given by
\begin{equation}
	u_*(x_1,x_2,x_3) = 212.5755x_1x_2+54.1296x_1x_3,
\end{equation}
which, in terms of the original state variables $(\theta,\dot\theta)$, is given by
\begin{equation}\label{PendController}
	u_*(\theta,\dot\theta) = 212.5755\cos(\theta)\sin(\theta) + 54.1296\cos(\theta)\dot\theta.
\end{equation}
We see in Figure~\ref{fig:PendControl} this control law leads to a quick jerk of the cart to the left and then right to swing the pendulum up, after which the cart ceases to move. Notice further that the controller \eqref{PendController} vanishes at both equilibria $(\theta,\dot\theta) = (0,0)$ and $(\theta,\dot\theta) =(\pi,0)$, meaning that the control is inactive at the hanging down ($\theta = \pi$) position. To circumvent this issue, a discontinuous control would be needed \cite{BHAT200063,mayhew_etal_2011_TAC,KALABIC2017293}, which is not explored herein. Nonetheless, this is little issue from a practical point of view as one may randomize an initial disturbance in the cart to throw oneself away from the hanging down state and then activate the control law to stabilize the system.


\section{Conclusion}
\label{sec:conclusion}

In this letter we have presented a simple, flexible, and efficient method for synthesizing control laws directly from data. We emphasize that our method does not require the intermediary step of identifying the nonlinear system dynamics or performing parameter estimation. Instead, we use the EDMD framework to produce an approximation of the Koopman operator on a dictionary of polynomial lifting functions. The result is a linear description of the dynamics in the lifted coordinates which can then be integrated with well-established CLF methods and convex SOS optimization routines. We therefore provide a two-step data-driven method that can be applied broadly to problems in control.   

Although we have explored the application to continuous-time data in this letter, the results can equivalently be applied to discrete-time processes as well by simply fixing $\tau = 1$. Moreover, the EDMD framework will equally approximate the Koopman operator for general stochastic processes with no modifications to the method and similar convergence guarantees \cite{bramburger2023auxiliary}. This means that our results could also be applied to data generated by stochastic systems, with the only minor difference being that the Lie derivative condition $\mathcal{L}V \leq 0$ is now considered in expectation \cite{kushner1967stochastic}. The ability to apply these methods to stochastic systems could be promising for overcoming the inevitable noise that is produced when gathering real-world data. A report on the applicability and noise-robustness of the method on laboratory data will be left to a follow-up investigation.

\bibliographystyle{IEEEtran}
\bibliography{data_control.bib}

\end{document}